\documentclass{notemat_arXiv}
\usepackage{amsmath}
\usepackage{amssymb}
\usepackage{nmmacro}

\usepackage{url}

\chardef\bslchar=`\\ 
\makeatletter
\newcommand{\addbslash}{\expandafter\@addbslash\string}
\def\@addbslash#1{\bslchar\@nobslash#1}
\newcommand{\nobslash}{\expandafter\@nobslash\string}
\def\@nobslash#1{\ifnum`#1=\bslchar\else#1\fi}
\newcommand{\ntt}{\normalfont\ttfamily}
\def\@boxorbreak{\leavevmode
  \ifmmode\hbox\else\ifdim\lastskip=\z@\penalty9999 \fi\fi}
\DeclareRobustCommand{\cs}[1]{\@boxorbreak{\ntt\addbslash#1\@empty}}
\makeatother
\begin{document}
\begin{article}
\begin{opening}
  \title{Representation growth and representation zeta functions of
    groups} 
  \author{Benjamin Klopsch\thanks{This survey reports on work which was
      partially supported by the following institutions: the Batsheva
      de Rothschild Fund, the EPSRC, the Mathematisches
      Forschungsinstitut Oberwolfach, the NSF and the Nuffield
      Foundation.}}
  \institute{Department of Mathematics, Royal Holloway University of London\\
    \email{Benjamin.Klopsch@rhul.ac.uk}}
  \runningauthor{B.\ Klopsch} \runningtitle{Representation growth and
    representation zeta functions of groups}
\begin{abstract}
  We give a short introduction to the subject of representation growth
  and representation zeta functions of groups, omitting all proofs.
  Our focus is on results which are relevant to the study of
  arithmetic groups in semisimple algebraic groups, such as the group
  $\textup{SL}_n(\mathbb{Z})$ consisting of $n \times n$ integer
  matrices of determinant $1$.  In the last two sections we state
  several results which were recently obtained in joint work with
  N.~Avni, U.~Onn and C.~Voll.
\end{abstract}

  \keywords{Representations, characters, arithmetic groups, $p$-adic
    Lie groups, zeta functions.}
  \classification{primary 22E55, secondary 22E50, 20F69}
\end{opening}


\section{Introduction} \label{sec:intro} Let $G$ be a group.  For $n
\in \mathbb{N}$, let $r_n(G)$ denote the number of isomorphism classes
of $n$-dimensional irreducible complex representations of~$G$.  We
suppose that $G$ is \emph{representation rigid}, i.e., that $r_n(G) <
\infty$ for all positive integers~$n$.

If the group $G$ is finite then $G$ is automatically representation
rigid and the sequence $r_n(G)$ has only finitely many non-zero terms,
capturing the distribution of irreducible character degrees of~$G$.
The study of finite groups by means of their irreducible character
degrees and conjugacy classes is a well established research area;
e.g., see \cite{Hu98} and references therein.  Interesting asymptotic
phenomena are known to occur when one considers the irreducible
character degrees of suitable infinite families of finite groups, for
instance, families of finite groups $H$ of Lie type as $\lvert H
\rvert$ tends to infinity; see~\cite{LiSh05}.

In the present survey we are primarily interested in the situation
where $G$ is infinite, albeit $G$ may sometimes arise as an inverse
limit of finite groups.  Two fundamental questions in this case are:
what are the arithmetic properties of the sequence $r_n(G)$, $n \in
\mathbb{N}$, and what is the asymptotic behaviour of $R_N(G) =
\sum_{n=1}^N r_n(G)$ as $N$ tends to infinity?  To a certain degree
this line of investigation is inspired by the subject of subgroup
growth and subgroup zeta functions which, in a similar way, is
concerned with the distribution of finite index subgroups; e.g., see
\cite{LuSe03,dSGr06}.

In order to streamline the investigation it is convenient to encode
the arithmetic sequence $r_n(G)$, $n \in \mathbb{N}$, in a suitable
generating function.  The \emph{representation zeta function} of $G$
is the Dirichlet generating function
\[
\zeta_G(s) = \sum_{n=1}^\infty r_n(G) \, n^{-s} \qquad (s \in
\mathbb{C}).
\] 
If the group $G$ is such that there is a one-to-one correspondence
between isomorphism classes of irreducible representations and
irreducible characters then, writing $\textup{Irr}(G)$ for the space
of irreducible characters of $G$, we can express the zeta function
also in the suggestive and slightly more algebraic form
\[
\zeta_G(s) = \sum_{\chi \in \textup{Irr}(G)} \chi(1)^{-s} \qquad (s
\in \mathbb{C}).
\]

The function $\zeta_G(s)$ is a suitable vehicle for studying the
distribution of character degrees of the group~$G$ whenever the
representation growth of $G$ is `not too fast', a condition which is
made precise in Section~\ref{sec:abscissa}.  Groups which meet this
requirement include, for instance, arithmetic groups in semisimple
algebraic groups with the Congruence Subgroup Property and open
compact subgroups of semisimple $p$-adic Lie groups.  In recent years,
several substantial results have been obtained concerning the
representation growth and representation zeta functions of these types
of groups; see~\cite{Ja06, LaLu08, AvKlOnVo10, Av11, AvKlOnVo12,
  AvKlOnVo12+, AvKlOnVo_arXiv, AvKlOnVo_preparation}.  In the present
survey we discuss some of these results and we indicate what kinds of
methods are involved in proving them.

\section{Finite groups of Lie type} Our primary focus is on infinite
groups, but it is instructive to briefly touch upon representation
zeta functions of finite groups of Lie type.  For example, the
representation theory of the general linear group
$\textup{GL}_2(\mathbb{F}_q)$ over a finite field $\mathbb{F}_q$ is
well understood and one deduces readily that
\begin{equation}\label{equ:GL2}
  \zeta_{\textup{GL}_2(\mathbb{F}_q)}(s) = (q-1) \left( 1 + q^{-s} +
    \tfrac{q-2}{2} (q+1)^{-s} + \tfrac{q}{2} (q-1)^{-s} \right).
\end{equation}
It is remarkable that the formula \eqref{equ:GL2} is uniform in $q$ in
the sense that both the irreducible character degrees and their
multiplicities can be expressed in terms of polynomials in $q$ over
the rational field~$\mathbb{Q}$.  In general, Deligne-Lusztig theory
provides powerful and sophisticated tools to study the irreducible
characters of finite groups of Lie type.  In~\cite{LiSh05}, Liebeck
and Shalev obtained, for instance, the following general asymptotic
result.

\begin{Theorem}[Liebeck and Shalev] \label{thm:finite-Lie-type} Let
  $L$ be a fixed Lie type and let $h$ be the Coxeter number of the
  corresponding simple algebraic group $\mathbf{G}$, i.e., $h + 1 =
  \dim \mathbf{G} / \textup{rk} \mathbf{G}$.  Then for the finite
  quasi-simple groups $L(q)$ of type $L$ over $\mathbb{F}_q$,
  \[
  \zeta_{L(q)}(s) \to
  \begin{cases}
      1 & \text{for $s \in \mathbb{R}_{> 2/h}$} \\
      \infty & \text{for $s \in \mathbb{R}_{< 2/h}$} 
    \end{cases}
    \qquad \text{as $q \to \infty$.}
    \]
\end{Theorem}

The Coxeter number $h$ is computed easily.  For example, for
$\mathbf{G} = \textup{SL}_n$ and $L(q) = \textup{SL}_n(\mathbb{F}_q)$
one has $h = n$.  In the smallest interesting case $n=2$ and for odd
$q$, the zeta function of $\textup{SL}_2(\mathbb{F}_q)$ is
\begin{equation}\label{equ:SL2}
  \zeta_{\textup{SL}_2(\mathbb{F}_q)}(s) = 1 + q^{-s} + \tfrac{q-3}{2}
  (q+1)^{-s} + \tfrac{q-1}{2} (q-1)^{-s}  + 2 (\tfrac{q+1}{2})^{-s} + 2
  (\tfrac{q-1}{2})^{-s},
\end{equation}
which is approximately the expression in~\eqref{equ:GL2} divided by
$(q-1)$.  From the explicit formula one can verify directly the
assertion of Theorem~\ref{thm:finite-Lie-type} in this special case.

\section{Abscissa of convergence and polynomial representation
  growth} \label{sec:abscissa} In Section~\ref{sec:intro} we
introduced the zeta function $\zeta_G(s)$ of a representation rigid
group $G$ as a formal Dirichlet series.  Clearly, if $G$ is finite --
or more generally if $G$ has only finitely many irreducible complex
representations -- then the Dirichlet polynomial $\zeta_G(s)$ defines
an analytic function on the entire complex plane.

Now suppose that $G$ is infinite and that $r_n(G)$ is non-zero for
infinitely many $n \in \mathbb{N}$.  Naturally, we are interested in
the convergence properties of $\zeta_G(s)$ for $s \in \mathbb{C}$.
The general theory of Dirichlet generating functions shows that the
region of convergence is always a right half plane of $\mathbb{C}$,
possibly empty, and that the resulting function is analytic.  If the
region of convergence is non-empty, one is also interested in
meromorphic continuation of the function to a larger part of the
complex plane.

The \emph{abscissa of convergence} $\alpha(G)$ of $\zeta_G(s)$ is the
infimum of all $\alpha \in \mathbb{R}$ such that the series
$\zeta_G(s)$ converges (to an analytic function) on the right half
plane $ \{ s \in \mathbb{C} \mid \textup{Re}(s) > \alpha \}$.  The
abscissa $\alpha(G)$ is finite if and only if $G$ has \emph{polynomial
  representation growth}, i.e., if $R_N(G) = \sum_{n=1}^N r_n(G)$
grows at most polynomially in~$N$.  In fact, if the growth sequence
$R_N(G)$, $N \in \mathbb{N}$, is unbounded then
\[
\alpha(G) = \limsup_{N \to \infty} \frac{\log R_N(G)}{\log N} 
\]
gives the \emph{polynomial degree of growth}: $R_N(G) =
O\left(N^{\alpha(G) + \varepsilon}\right)$ for every $\varepsilon >
0$.  

Two fundamental problems in the subject are: to characterise groups of
polynomial representation growth -- motivated by Gromow's celebrated
theorem on groups of polynomial word growth -- and to link the actual
value of the abscissa of convergence $\alpha(G)$ of a group $G$ to
structural properties of $G$.  In general these questions are still
very much open.  However, in the context of semisimple algebraic
groups and their arithmetic subgroups a range of results have been
obtained.  A selection of these are discussed in the following
sections.

\section{Witten zeta functions} \label{sec:Witten} In \cite{Wi91},
Witten initiated in the context of quantum gauge theories the study of
certain representation zeta functions.  Let $\mathbf{G}$ be a
connected, simply connected, complex almost simple algebraic group and
let $G = \mathbf{G}(\mathbb{C})$.  It is natural to focus on rational
representations of the algebraic group $G$ and one can show that $G$
is representation rigid in this restricted sense.  The \emph{Witten
  zeta function} $\zeta_G(s)$ counts irreducible rational
representations of the complex algebraic group $G$.  These zeta
functions also appear naturally as archimedean factors of
representation zeta functions of arithmetic groups, as explained in
Section~\ref{sec:lattices}.

For example, the group $\textup{SL}_2(\mathbb{C})$ has a unique
irreducible representation of each possible degree.  Hence
\[
\zeta_{\textup{SL}_2(\mathbb{C})}(s) = \sum_{n=1}^\infty n^{-s},
\]
the famous Riemann zeta function.  In particular, the abscissa
of convergence is~$1$ and there is a meromorphic continuation to the
entire complex plane.

In general, the irreducible representations $V_\lambda$ of $G$ are
parametrised by their highest weights $\lambda = \sum_{i=1}^r a_i
\overline{\omega}_i$, where $\overline{\omega}_1$, \ldots,
$\overline{\omega}_r$ denote the fundamental weights and the
coefficients $a_1, \dots, a_r$ range over all non-negative integers.
Moreover, $\dim V_\lambda$ is given by the Weyl dimension formula.
By a careful analysis, Larsen and Lubotzky prove in \cite{LaLu08} the
following result.

\begin{Theorem}[Larsen and Lubotzky] \label{thm:Witten}
  Let $\mathbf{G}$ be a connected, simply connected, complex almost
  simple algebraic group and let $G = \mathbf{G}(\mathbb{C})$.  Then
  $\alpha(G) = 2/h$, where $h$ is the Coxeter number of $\mathbf{G}$.
\end{Theorem}

It is known that Witten zeta functions can be continued
meromorphically to the entire complex plane.  Further analytic
properties of these functions, such as the location of singularities
and functional relations, have been investigated in some detail using
multiple zeta functions; e.g., see~\cite{Ma03,KoMaTs10}.  It is
remarkable that the same invariant $2/h$ features in
Theorems~\ref{thm:finite-Lie-type} and \ref{thm:Witten}.  Currently
there appears to be no conceptual explanation for this.
 
\section{The group $\textup{SL}_2(R)$ for discrete valuation rings
  $R$} \label{sec:SL2} If $G$ is a topological group it is natural to
focus attention on continuous representations.  A finitely generated
profinite group $G$ is representation rigid in this restricted sense
if and only if it is \emph{FAb}, i.e., if every open subgroup $H$ of
$G$ has finite abelianisation~$H/[H,H]$.  This is a consequence of
Jordan's theorem on abelian normal subgroups of bounded index in
finite linear groups.  We tacitly agree that the representation zeta
function $\zeta_G(s)$ of a finitely generated FAb profinite group $G$
counts irreducible continuous complex representations of~$G$.

Let $R$ be a complete discrete valuation ring, with residue field
$\mathbb{F}_q$ of odd characteristic.  This means that $R$ is either a
finite integral extension of the ring of $p$-adic integers
$\mathbb{Z}_p$ for some prime $p$ or a formal power series ring
$\mathbb{F}_q[\![t]\!]$ over a finite field of cardinality~$q$.

In~\cite{Ja06}, Jaikin-Zapirain showed by a hands-on computation of
character degrees that the representation zeta function
$\zeta_{\textup{SL}_2(R)}(s)$ equals
\[
\zeta_{\textup{SL}_2(\mathbb{F}_q)}(s) + \left( 4q \left(
    \tfrac{q^2-1}{2} \right)^{-s} + \tfrac{q^2-1}{2} (q^2-q)^{-s} +
  \tfrac{(q-1)^2}{2} (q^2 + q)^{-s} \right)/(1-q^{1-s}),
\]
where the Dirichlet polynomial
$\zeta_{\textup{SL}_2(\mathbb{F}_q)}(s)$ is described in
\eqref{equ:SL2}.  It is remarkable that the above formula is uniform
in $q$, irrespective of the characteristic, absolute ramification
index or isomorphism type of the ring~$R$.  In the case where $R$ has
characteristic $0$, Lie-theoretic techniques combined with Clifford
theory can be used to gain an insight into the features of this
specific example which hold more generally; see Sections~\ref{sec:FAb}
and \ref{sec:AKOV}.

Clearly, the explicit formula for the function
$\zeta_{\textup{SL}_2(R)}(s)$ provides a meromorphic extension to the
entire complex plane.  The abscissa of convergence is $1$ and, in view
of Theorems~\ref{thm:finite-Lie-type} and \ref{thm:Witten}, this value
could be interpreted as $2/h$, the Coxeter number of $\textup{SL}_2$
being $h=2$.  But such an interpretation is misleading, as can be seen
from the following general result obtained in~\cite{LaLu08}.

\begin{Theorem}[Larsen and Lubotzky] \label{thm:isotropic} Let
  $\mathbf{G}$ be a simple algebraic group over a non-archimedean
  local field~$F$.  Suppose that $\mathbf{G}$ is $F$-isotropic, i.e.,
  $\textup{rk}_F \mathbf{G} \geq 1$.  Let $H$ be a compact open
  subgroup of $\mathbf{G}(F)$.  Then $\alpha(H) \geq 1/15$.
\end{Theorem}

Taking $\mathbf{G} = \textup{SL}_n$ and $F = \mathbb{Q}_p$, we may
consider the compact $p$-adic Lie groups
$\textup{SL}_n(\mathbb{Z}_p)$. For these groups $2/h = 2/n \to 0$ as
$n \to \infty$, whereas $\alpha(\textup{SL}_n(\mathbb{Z}_p))$ is
uniformly bounded away from $0$.  Currently, the only explicit values
known for $\alpha(\textup{SL}_n(\mathbb{Z}_p))$ are: $1$ for $n=2$ (as
seen above), and $2/3$ for $n=3$ (see \cite{AvKlOnVo12+}).
Unfortunately, these do not yet indicate the general behaviour.

\section{FAb compact $p$-adic Lie groups} \label{sec:FAb} Let $G$ be a
compact $p$-adic Lie group.  Then one associates to $G$ a
$\mathbb{Q}_p$-Lie algebra as follows.  The group $G$ contains a
uniformly powerful open pro-$p$ subgroup $U$.  By the theory of
powerful pro-$p$ groups, $U$ gives rise to a $\mathbb{Z}_p$-Lie
lattice $L = \log(U)$ and the induced $\mathbb{Q}_p$-Lie algebra
$\mathcal{L}(G) = \mathbb{Q}_p \otimes_{\mathbb{Z}_p} L$ does not
depend on the specific choice of $U$.  It is a fact that $G$ is FAb if
and only if $\mathcal{L}(G)$ is perfect, i.e., if
$[\mathcal{L}(G),\mathcal{L}(G)] = \mathcal{L}(G)$.  Conversely, for
any $\mathbb{Q}_p$-Lie algebra $\mathcal{L}$ one can easily produce
compact $p$-adic Lie groups $G$ such that $\mathcal{L}(G) =
\mathcal{L}$, using the exponential map.  This supplies a large class
of compact $p$-adic Lie groups which are FAb and hence have polynomial
representation growth.

Using the Kirillov orbit method and techniques from model theory,
Jaikin-Zapirain established in \cite{Ja06} that the representation
zeta function of a FAb compact $p$-adic analytic pro-$p$ group can
always be expressed as a rational function in $p^{-s}$ over
$\mathbb{Q}$.  More generally, he proved the following result, which
is illustrated by the explicit example $G = \textup{SL}_2(R)$ given in
Section~\ref{sec:SL2}.

\begin{Theorem}[Jaikin-Zapirain]
  Let $G$ be an FAb compact $p$-adic Lie group, and suppose that
  $p>2$.  Then there are finitely many positive integers $n_1, \ldots,
  n_k$ and rational functions $f_1, \ldots, f_k \in \mathbb{Q}(X)$
  such that
  \[
  \zeta_G(s) = \sum_{i=1}^k f_i(p^{-s}) \, n_i^{-s}.
  \]
\end{Theorem}

In particular, the theorem shows that the zeta function of a FAb
compact $p$-adic Lie group $G$ extends meromorphically to the entire
complex plane.  The invariant $\alpha(G)$ is the largest real part of
a pole of $\zeta_G(s)$.  It is natural to investigate the whole
spectrum of poles and zeros of $\zeta_G(s)$.

Currently, very little is known about the location of the zeros of
representation zeta functions.  In 2010 Kurokawa and Kurokawa observed
from the explicit formula given in Section~\ref{sec:SL2} that
$\zeta_{\textup{SL}_2(\mathbb{Z}_p)}(s) = 0$ for $s \in \{-1,-2\}$.
We note that if $G$ is a finite group then $\zeta_G(-2) = \sum_{\chi
  \in \textup{Irr}(G)} \chi(1)^2 = \lvert G \rvert$.  Based on this
fact and the results in \cite{Ja06} one can prove the following
general result.

\begin{Theorem}[Jaikin-Zapirain and Klopsch]
  Let $G$ be an infinite FAb compact $p$-adic Lie group and suppose
  that $p > 2$.  Then $\zeta_G(-2) = 0$.
\end{Theorem}

\section{Rational representations of the infinite cyclic group}

Before considering arithmetic subgroups of semisimple algebraic
groups, let us look at representations of the simplest infinite group,
i.e., the infinite cyclic group $C_\infty$.  The group $C_\infty$ has
already infinitely many $1$-dimensional representations.  Hence in
order to say anything meaningful we need to slightly adapt our basic
definitions.

We make two modifications: firstly let us only consider
representations with finite image and secondly let us consider
irreducible representations over $\mathbb{Q}$ rather than
$\mathbb{C}$.  More precisely, for any finitely generated nilpotent
group $\Gamma$ let $\hat r_n^\mathbb{Q}(\Gamma)$ denote the number of
$n$-dimensional irreducible representations of $\Gamma$ over
$\mathbb{Q}$ with finite image.  Then it turns out that $\hat
r_n^\mathbb{Q}(\Gamma)$ is finite for every $n \in \mathbb{N}$ and we
can define the \emph{$\mathbb{Q}$-rational representation zeta
  function}
\[
\zeta^\mathbb{Q}_\Gamma(s) = \sum_{n=1}^\infty \hat r_n^\mathbb{Q}(\Gamma)
\, n^{-s}.
\]
Using that every finite nilpotent group is the direct product of its
Sylow $p$-subgroups and basic facts from character theory, one can show
that $\zeta^\mathbb{Q}_\Gamma(s)$ admits an \emph{Euler product
  decomposition}
\begin{equation} \label{equ:Euler} \zeta^\mathbb{Q}_\Gamma(s) =
  \prod_{p \text{ prime}} \zeta^\mathbb{Q}_{\Gamma,p}(s),
\end{equation}
where for each prime $p$ the local factor
$\zeta^\mathbb{Q}_{\Gamma,p}(s) = \sum_{k=0}^\infty \hat
r_{p^k}^\mathbb{Q}(\Gamma) \, p^{-ks}$, enumerating irreducible
representations of $p$-power dimension, can be re-interpreted as the
$\mathbb{Q}$\nobreakdash-rational representation zeta function of the
pro-$p$ completion $\widehat{\Gamma}_p$ of $\Gamma$.  For more details
and deeper results in this direction we refer to the forthcoming
article~\cite{GoJaKl_preparation}.

Let us now return to the simplest case: $\Gamma = C_\infty$, the
infinite cyclic group.  Since the group $C_\infty$ is abelian, its
irreducible representations over $\mathbb{Q}$ with finite image can be
effectively described by means of Galois orbits of irreducible complex
characters.  In the general setting, one would also need to keep track
of Schur indices featuring in the computation of
$\zeta^\mathbb{Q}_{\Gamma,2}(s)$.  A short analysis yields
\[
\zeta_{C_\infty}^\mathbb{Q}(s) = \sum_{m=1}^\infty \varphi(m)^{-s},
\]
where $\varphi$ denotes Euler's function familiar from elementary
number theory.

The Dirichlet series $\psi(s) = \sum_{m=1}^\infty \varphi(m)^{-s}$ is
of independent interest in analytic number theory and has been studied
by many authors; e.g., see \cite{Ba72}.  The Euler product
decomposition \eqref{equ:Euler} can be established directly
\[
\psi(s) = \prod_\text{$p$ prime} \left( 1 + (p-1)^{-s} / (1- p^{-s})
\right).
\]
The abscissa of convergence of $\psi(s)$, which can be interpreted as
the degree $\alpha^\mathbb{Q}(C_\infty)$ of $\mathbb{Q}$-rational
representation growth, is equal to $1$.  In fact, writing
\[ 
\psi(s) = \underbrace{\prod_\text{$p$ prime} \left( 1 + (p-1)^{-s} -
    p^{-s} \right)}_{\text{converges for $\textup{Re}(s) > 0$}} \cdot
\underbrace{\prod_\text{$p$ prime} (1- p^{-s})^{-1}}_{\text{Riemann
    zeta function $\zeta(s)$}},
\]
one sees that $\psi(s)$ admits a meromorphic continuation to
$\textup{Re}(s) > 0$ (but not to the entire complex plane) and has a
simple pole at $s=1$ with residue $c = \zeta(2) \zeta(3) / \zeta(6) =
1.9435964\ldots$ This yields very precise asymptotics for the
$\mathbb{Q}$-rational representation growth of $C_\infty$; in
particular,
\[
\sum_{n=1}^N \hat r_n^\mathbb{Q}(C_\infty) = \# \{ m \mid \varphi(m)
\leq N \} \sim cN \qquad \text{as $N \to \infty$}.
\]
One may regard this simple case and its beautiful connections to
classical analytic number theory as a further motivation for studying
representation zeta functions of arithmetic groups.

\section{Arithmetic lattices in semisimple
  groups} \label{sec:lattices} In this section we turn our attention
to lattices in semisimple locally compact groups.  These lattices are
discrete subgroups of finite co-volume and often, but not always, have
arithmetic origin.  For instance, $\textup{SL}_n(\mathbb{Z})$ is an
arithmetic lattice in the real Lie group $\textup{SL}_n(\mathbb{R})$.
More generally, let $\Gamma$ be an arithmetic irreducible lattice in a
semisimple locally compact group $G$ of characteristic $0$.  Then
$\Gamma$ is commensurable to $\mathbf{G}({\scriptstyle
  \mathcal{O}}_S)$, where $\mathbf{G}$ is a connected, simply
connected absolutely almost simple algebraic group defined over a
number field $k$ and ${\scriptstyle \mathcal{O}}_S$ is the ring of
$S$-integers for a finite set $S$ of places of~$k$.  By a theorem
going back to Borel and Harish-Chandra, any such
$\mathbf{G}({\scriptstyle \mathcal{O}}_S)$ forms an irreducible
lattice in the semisimple locally compact group $G = \prod_{\wp \in S}
\mathbf{G}(k_\wp)$ under the diagonal embedding, as long as $S$ is
non-empty and contains all archimedean places $\wp$ such that
$\mathbf{G}(k_\wp)$ is non-compact.  Examples of this construction are
$\textup{SL}_n(\mathbb{Z}[\sqrt{2}]) \subseteq
\textup{SL}_n(\mathbb{R}) \times \textup{SL}_n(\mathbb{R})$ and
$\textup{SL}_n(\mathbb{Z}[1/p]) \subseteq \textup{SL}_n(\mathbb{R})
\times \textup{SL}_n(\mathbb{Q}_p)$.  Margulis has shown that in the
higher rank situation all irreducible lattices are arithmetic and
arise in this way.  For precise notions and a more complete
description see~\cite{Ma91}.

Throughout the following we assume, for simplicity of notation, that
$\Gamma = \mathbf{G}({\scriptstyle \mathcal{O}}_S)$ as above.  In
\cite{LuMa04}, Lubotzky and Martin gave a characterisation of
arithmetic groups of polynomial representation growth, linking them to
the classical Congruence Subgroup Problem.

\begin{Theorem}[Lubotzky and Martin]
  Let $\Gamma$ be an arithmetic group as above.  Then $\alpha(\Gamma)$
  is finite if and only if $\Gamma$ has the Congruence Subgroup
  Property.
\end{Theorem}

The group $\Gamma$ has the Congruence Subgroup Property (CSP) if,
essentially, all its finite index subgroups arise from the arithmetic
structure of the group.  Technically, this means that the congruence
kernel $\ker (\widehat{\mathbf{G}({\scriptstyle \mathcal{O}}_S)}
\rightarrow \overline{\mathbf{G}({\scriptstyle \mathcal{O}}_S)} )$ is
finite; here $\widehat{\mathbf{G}({\scriptstyle \mathcal{O}}_S)}$ is
the profinite completion and $\overline{\mathbf{G}({\scriptstyle
    \mathcal{O}}_S)} \cong \prod_{\mathfrak{p} \not \in S}
\mathbf{G}({\scriptstyle \mathcal{O}}_\mathfrak{p})$, with
$\mathfrak{p}$ running over non-archimedean places, denotes the
congruence completion of $\mathbf{G}({\scriptstyle \mathcal{O}}_S)$.
For instance, it was shown by Bass-Lazard-Serre and Mennicke that the
group $\textup{SL}_n(\mathbb{Z})$ has the CSP if and only if $n \geq
3$.  That $\textup{SL}_2(\mathbb{Z})$ does not have the CSP was
discovered by Fricke and Klein.  Retrospectively this is not
surprising, because $\textup{SL}_2(\mathbb{Z})$ contains a free
subgroup of finite index.  We refer to \cite{PrRa10} for a
comprehensive survey of the Congruence Subgroup Problem, i.e., the
problem to decide precisely which arithmetic groups have the CSP.

Suppose that $\Gamma$ has the CSP.  Using Margulis' super-rigidity
theorem, Larsen and Lubotzky derived in \cite{LaLu08} an Euler product
decomposition for $\zeta_\Gamma(s)$, which takes a particularly simple
form whenever the congruence kernel is trivial.

\begin{Theorem}[Larsen and Lubotzky]
  Let $\Gamma$ be an arithmetic group as above and suppose that
  $\Gamma$ has the CSP.  Then $\zeta_\Gamma(s)$ admits an Euler
  product decomposition.  In particular, if the congruence kernel for
  $\Gamma = \mathbf{G}({\scriptstyle \mathcal{O}}_S)$ is trivial then
  \begin{equation}\label{equ:EulerLaLu}
    \zeta_{\Gamma}(s) =
    \zeta_{\mathbf{G}(\mathbb{C})}(s)^{[k:\mathbb{Q}]} \,
    \prod_\text{$\mathfrak{p} \not \in S$} \zeta_{\mathbf{G}({\scriptscriptstyle
        \mathcal{O}}_\mathfrak{p})}(s).
  \end{equation}
\end{Theorem}

For instance, for the groups $\textup{SL}_n(\mathbb{Z})$, $n \geq 3$,
the Euler product takes the form
\[
\zeta_{\textup{SL}_n(\mathbb{Z})}(s) =
\zeta_{\textup{SL}_n(\mathbb{C})}(s) \, \prod_\text{$p$ prime}
\zeta_{\textup{SL}_n(\mathbb{Z}_p)}(s).
\]
In Sections~\ref{sec:Witten} and \ref{sec:SL2} we already encountered
individually the factors of these Euler products:
$\zeta_{\mathbf{G}(\mathbb{C})}(s)$ is the Witten zeta function
capturing rational representations of the algebraic group
$\mathbf{G}(\mathbb{C})$ and, for each $\mathfrak{p}$, the function
$\zeta_{\mathbf{G}({\scriptscriptstyle \mathcal{O}}_\mathfrak{p})}(s)$
enumerates continuous representations of the compact $p$-adic Lie
group $\mathbf{G}({\scriptscriptstyle \mathcal{O}}_\mathfrak{p})$.
Larsen and Lubotzky's results for the abscissae of convergence of
these local zeta functions include Theorems~\ref{thm:Witten} and
\ref{thm:isotropic} stated above.

Regarding the abscissa of convergence of the global representation
zeta function, Avni employed in~\cite{Av11} model-theoretic techniques
to prove that the abscissa of convergence of $\zeta_\Gamma(s)$ is
always a rational number.  In \cite{LaLu08}, Larsen and Lubotzky made
the following conjecture, which can be regarded as a refinement of
Serre's conjecture on the Congruence Subgroup Problem.

\begin{Conjecture}[Larsen and Lubotzky]
  Let $G$ be a higher-rank semisimple locally compact group. Then, for
  any two irreducible lattices $\Gamma_1$ and $\Gamma_2$ in $G$,
  $\alpha(\Gamma_1) = \alpha(\Gamma_2)$.
\end{Conjecture}

Roughly speaking, the conjecture states that the ambient semisimple
locally compact group does not only control whether lattices contained
in it have the CSP (as in Serre's conjecture), but also what their
polynomial degree of representation growth is.  A concrete example of
a lattice in $\textup{SL}_n(\mathbb{R})$ which is rather different
from the most familiar one $\textup{SL}_n(\mathbb{Z})$ is the special
unitary group $\textup{SU}_n(\mathbb{Z}[\sqrt{2}],\mathbb{Z})$,
consisting of all matrices $A = (a_{ij})$ over the ring
$\mathbb{Z}[\sqrt{2}]$ with $\det A = 1$ and $A^{-1} =
(a_{ji}^\sigma)$, where $\sigma$ is the Galois automorphism of
$\mathbb{Q}(\sqrt{2})$ swapping $\sqrt{2}$ and $-\sqrt{2}$.
 
\section{New results for arithmetic groups and compact
  $p$\nobreakdash-adic Lie groups} \label{sec:AKOV} The short
announcement~\cite{AvKlOnVo10} summarises a number of results obtained
recently by the author in joint work with Avni, Onn and Voll.  Details
are appearing in \cite{AvKlOnVo12,AvKlOnVo12+,AvKlOnVo_preparation}.
The toolbox which we use to prove our results comprises a variety of
techniques which can only be hinted at: they include, for instance,
the Kirillov orbit method for $p$-adic analytic pro-$p$ groups,
methods from $\mathfrak{p}$-adic integration and the study of
generalised Igusa zeta functions, the theory of sheets of simple Lie
algebras, resolution of singularities in characteristic $0$, aspects
of the Weil conjectures regarding zeta functions of smooth projective
varieties over finite fields, approximative and exact Clifford theory.

In summary our main results are
\begin{itemize}
\item a global \emph{Denef formula} for the zeta functions of
  principal congruence subgroups of compact $p$-adic Lie groups, such
  as $\textup{SL}_n^m(\mathbb{Z}_p) \subseteq
  \textup{SL}_n(\mathbb{Z}_p)$;
\item local \emph{functional equations} for the zeta functions of
  principal congruence subgroups of compact $p$-adic Lie groups, such
  as $\textup{SL}_n^m(\mathbb{Z}_p) \subseteq
  \textup{SL}_n(\mathbb{Z}_p)$;
\item \emph{candidate pole sets} for the non-archimedean factors
  occurring in the Euler product \eqref{equ:EulerLaLu}, e.g., the zeta
  functions $\zeta_{\textup{SL}_n(\mathbb{Z}_p)}(s)$;
\item \emph{explicit formulae} for the zeta functions of compact
  $p$-adic Lie groups of type $A_2$, such as
  $\textup{SL}_3(\mathbb{Z}_p)$ and
  $\textup{SU}_3(\mathfrak{O},\mathbb{Z}_p)$ for unramified
  $\mathfrak{O}$;
\item \emph{meromorphic continuation} of zeta functions and a precise
  \emph{asymptotic description} of the representation growth for
  arithmetic groups of type $A_2$, such as
  $\textup{SL}_3(\mathbb{Z})$.
\end{itemize}

These results are clearly relevant in the context of the Euler
product~\eqref{equ:EulerLaLu}.  Moreover, a large part of our work
applies in a more general context than discussed so far.  We recall
from Section~\ref{sec:FAb} that a compact $p$-adic Lie group $G$ is
representation rigid if and only if its $\mathbb{Q}_p$-Lie algebra
$\mathcal{L}(G)$ is perfect.  Let $k$ be a number field, and let
${\scriptstyle \mathcal{O}}$ be its ring of integers.  Let $\Lambda$
be an ${\scriptstyle \mathcal{O}}$-Lie lattice such that $k
\otimes_{\scriptstyle \mathcal{O}} \Lambda$ is perfect of
dimension~$d$.  Let $\mathfrak{o}$ be the completion ${\scriptstyle
  \mathcal{O}}_\mathfrak{p}$ of ${\scriptstyle \mathcal{O}}$ at a
non-archimedean place $\mathfrak{p}$.  Let $\mathfrak{O}$ be a finite
integral extension of $\mathfrak{o}$, corresponding to a place
$\mathfrak{P}$ lying above~$\mathfrak{p}$.  For $m \in \mathbb{N}$,
let $\mathfrak{g}^m(\mathfrak{O})$ denote the $m$th principal
congruence Lie sublattice of the $\mathfrak{O}$-Lie lattice
$\mathfrak{O} \otimes_{\scriptstyle \mathcal{O}} \Lambda$.  For
sufficiently large $m$, let $\mathsf{G}^m(\mathfrak{O})$ be the
$p$-adic analytic pro-$p$ group $\exp(\mathfrak{g}^m(\mathfrak{O}))$.

Using the Kirillov orbit method for permissible
$\mathsf{G}^m(\mathfrak{O})$, e.g., $\textup{SL}_n^1(\mathbb{Z}_p)$,
we can `linearise' the problem of enumerating irreducible characters
of the group $\mathsf{G}^m(\mathfrak{O})$ by their degrees.  We then
set up a generalised Igusa zeta function, i.e., a $p$-adic integral of
the form
\[
\mathcal{Z}_{\mathfrak{O}}(r,t) = \int_{(x,\mathbf{y}) \in
  V(\mathfrak{O})} \lvert x \rvert_\mathfrak{P}^t \prod_{j=1}^{\lfloor
  d/2 \rfloor} \frac{\lVert F_j(\mathbf{y}) \cup
  F_{j-1}(\mathbf{y})x^2 \rVert_\mathfrak{P}^r}{\lVert
  F_{j-1}(\mathbf{y}) \rVert_\mathfrak{P}^r} \, d\mu(x,\mathbf{y}),
\]
where $V(\mathfrak{O})\subset \mathfrak{O}^{d+1}$ is a union of cosets
modulo $\mathfrak{P}$, $F_j(\mathbf{Y}) \subset {\scriptstyle
  \mathcal{O}}[\mathbf{Y}]$ are polynomial sets defined in terms of
the structure constants of the underlying ${\scriptstyle
  \mathcal{O}}$-Lie lattice~$\Lambda$, $\|\cdot \|_\mathfrak{P}$ is
the $\mathfrak{P}$-adic maximum norm and $\mu$ is the additive Haar
measure on $\mathfrak{O}^{d+1}$ with $\mu(\mathfrak{O}^{d+1})=1$.  The
integral $\mathcal{Z}_{\mathfrak{O}}(r,t)$ allows us to treat
`uniformly' the representation zeta functions of the different groups
$\exp(\mathsf{G}^m(\mathfrak{O}))$ arising from the global
${\scriptstyle \mathcal{O}}$-Lie lattice $\Lambda$ under variation of
the place $\mathfrak{p}$ of ${\scriptstyle \mathcal{O}}$, the local
ring extension $\mathfrak{O}$ of ${\scriptstyle
  \mathcal{O}}_\mathfrak{p}$ and the congruence level $m$.  In
particular, we derive from our analysis a Denef formula and local
functional equations.

\begin{Theorem}[Avni, Klopsch, Onn and Voll~\cite{AvKlOnVo12+}]
  In the setup described, there exist $r \in \mathbb{N}$ and a
  rational function $R(X_1,\dots,X_r,Y)\in\mathbb{Q}(X_1,\dots,X_r,Y)$
  such that for almost every non-archimedean place $\mathfrak{p}$ of $k$ the
  following holds.

  There are algebraic integers $\lambda_1, \dots, \lambda_r$ such that
  for all finite extensions $\mathfrak{O}$ of $\mathfrak{o} =
  {\scriptstyle \mathcal{O}}_\mathfrak{p}$ and all permissible $m$ one has
  \[
  \zeta_{\mathsf{G}^m(\mathfrak{O})}(s) = q_\mathfrak{p}^{f dm}R(\lambda_1^f,
  \dots, \lambda_r^f,q_\mathfrak{p}^{-fs}),
  \]
  where $q_\mathfrak{p}$ is the residue field cardinality of
  $\mathfrak{o}$, $f$ denotes the inertia degree of $\mathfrak{O}$
  over $\mathfrak{o}$ and $d = \dim_k(k \otimes_{\scriptscriptstyle
    \mathcal{O}} \Lambda)$.  Moreover, there is the functional
  equation
  \[
  \zeta_{\mathsf{G}^m(\mathfrak{O})}(s)|_{\substack{q_\mathfrak{p}\rightarrow
      q_\mathfrak{p}^{-1} \\ \lambda_i\rightarrow
      \lambda_i^{-1}}}=q_\mathfrak{p}^{fd(1-2m)}\zeta_{\mathsf{G}^m(\mathfrak{O})}(s).
  \]
\end{Theorem}

Furthermore, we obtain candidate pole sets and we show that, locally,
abscissae of convergence are monotone under ring extensions.

\begin{Theorem}[Avni, Klopsch, Onn and Voll~\cite{AvKlOnVo12+}]
  In the setup described, there exists a finite set $P \subset
  \mathbb{Q}_{>0}$ such that the following is true.
  
  For all non-archimedean places $\mathfrak{p}$ of $k$, all finite extensions
  $\mathfrak{O}$ of $\mathfrak{o} = {\scriptstyle \mathcal{O}}_\mathfrak{p}$ and
  all permissible $m$ one has
  \[
  \left\{ \textup{Re}(z) \mid z \in \mathbb{C} \text{ a pole of
    }\zeta_{\mathsf{G}^m(\mathfrak{O})}(s) \right\} \subseteq P.
  \]
  In particular, one has $\alpha(\mathsf{G}^m(\mathfrak{O})) \leq \max
  P$, and equality holds for a set of positive Dirichlet density.

  Furthermore, if $\mathfrak{p}$ is any non-archimedean place of $k$ and if
  ${\scriptstyle \mathcal{O}}_\mathfrak{p} = \mathfrak{o} \subseteq
  \mathfrak{O}_1 \subseteq \mathfrak{O}_2$ is a tower of finite ring
  extensions, then for every permissible $m$ one has
  \[
  \alpha(\mathsf{G}^m(\mathfrak{O}_1)) \leq
  \alpha(\mathsf{G}^m(\mathfrak{O}_2)).
  \]
\end{Theorem}

By a more detailed study of groups of type $A_2$, we obtain the
following theorems addressing, in particular, the conjecture of Larsen
and Lubotzky stated in Section~\ref{sec:lattices}.  Analysing the
unique subregular sheet of the Lie algebra
$\mathfrak{sl}_3(\mathbb{C})$ and using approximative Clifford theory,
we prove the next result.

\begin{Theorem}[Avni, Klopsch, Onn and
  Voll~\cite{AvKlOnVo12+}] \label{thm:A2} Let $\Gamma$ be an
  arithmetic subgroup of a connected, simply connected simple
  algebraic group of type $A_2$ defined over a number field.  If
  $\Gamma$ has the CSP, then $\alpha(\Gamma) = 1$.
\end{Theorem}
 
Employing exact Clifford theory, we obtain the following more detailed
result for the special linear group $\textup{SL}_3({\scriptstyle
  \mathcal{O}})$ over the ring of integers of a number field.

\begin{Theorem}[Avni, Klopsch, Onn and Voll~\cite{AvKlOnVo_preparation}]
  Let ${\scriptstyle \mathcal{O}}$ be the ring of integers of a number
  field $k$.  Then there exists $\varepsilon > 0$ such that the
  representation zeta function of $\textup{SL}_3({\scriptstyle
    \mathcal{O}})$ admits a meromorphic continuation to the half-plane
  $\{ s \in \mathbb{C} \mid \textup{Re}(s) > 1 - \varepsilon \}$.  The
  continued function is analytic on the line $\{s \in \mathbb{C} \mid
  \textup{Re}(s) = 1 \}$, except for a double pole at $s = 1$.

  Consequently, there is a constant $c \in \mathbb{R}_{>0}$ such that
  \[
  R_N(\textup{SL}_3({\scriptstyle \mathcal{O}})) = \sum_{n=1}^N
  r_n(\textup{SL}_3({\scriptstyle \mathcal{O}})) \sim c \cdot N (log
  N) \qquad \text{as $N \to \infty$.}
  \]
\end{Theorem}

A key step in proving this result consists in deriving explicit
formulae for the representation zeta function of groups
$\textup{SL}_3(\mathfrak{o})$, where $\mathfrak{o}$ is a compact
discrete valuation ring of characteristic $0$ and residue field
characteristic different from~$3$.  In fact, we also derive similar
results for special unitary groups
$\textup{SU}_3(\mathcal{O},{\scriptstyle \mathcal{O}})$.

\section{New results regarding the conjecture of Larsen and
  Lubotzky}

Very recently, in joint work with Avni, Onn and Voll we prove the
following theorem in connection with the conjecture of Larsen and
Lubotzky which is stated in Section~\ref{sec:lattices}.

\begin{Theorem}[Avni, Klopsch, Onn and
  Voll~\cite{AvKlOnVo_arXiv}] \label{thm:phi} Let $\Phi$ be an
  irreducible root system.  Then there exists a constant $\alpha_\Phi$
  such that for every number field $k$ with ring of integers
  ${\scriptstyle \mathcal{O}}$, every finite set $S$ of places of $k$
  and every connected, simply connected absolutely almost simple
  algebraic group $\mathbf{G}$ over $k$ with absolute root system
  $\Phi$ the following holds.

  If the arithmetic group $\mathbf{G}({\scriptstyle \mathcal{O}}_S)$
  has polynomial representation growth, then
  $\alpha(\mathbf{G}({\scriptstyle \mathcal{O}}_S)) = \alpha_\Phi$.
\end{Theorem}

On the one hand, Theorem~\ref{thm:phi} is weaker than the conjecture
of Larsen and Lubotzky, because it does not resolve Serre's conjecture
on the Congruence Subgroup Problem.  However, Serre's conjecture is
known to be true in many cases and we have the following corollary.

\begin{Corollary}
  Serre's conjecture on the Congruence Subgroup Problem implies Larsen
  and Lubotzky's conjecture on the degrees of representation growth of
  lattices in higher rank semisimple locally compact groups.
\end{Corollary}

On the other hand, Theorem~\ref{thm:phi} is stronger than the
conjecture of Larsen and Lubotzky, because it shows that many
arithmetic groups with the CSP have the same degree of representation
growth, even when they do not embed as lattices into the same
semisimple locally compact group.  For instance, fixing $\Phi$ of type
$A_{n-1}$ for some $n \geq 3$, all of the following groups (for which
we also display their embeddings as lattices into semisimple locally
compact groups) have the same degree of representation growth:
\begin{enumerate}
\item $\textup{SL}_n(\mathbb{Z}) \subseteq \textup{SL}_n(\mathbb{R})$,
\item $\textup{SL}_n(\mathbb{Z}[\sqrt{2}]) \subseteq
  \textup{SL}_n(\mathbb{R}) \times \textup{SL}_n(\mathbb{R})$,
\item $\textup{SL}_n(\mathbb{Z}[i]) \subseteq
  \textup{SL}_n(\mathbb{C})$,
\item $\textup{SL}_n(\mathbb{Z}[1/p]) \subseteq
  \textup{SL}_n(\mathbb{R}) \times \textup{SL}_n(\mathbb{Q}_p)$,
\item $\textup{SU}_n(\mathbb{Z}[\sqrt{2}],\mathbb{Z}) \subseteq
  \textup{SL}_n(\mathbb{R})$.
\end{enumerate}

Presently, the only known explicit values of $\alpha_\Phi$ are: $2$
for $\Phi$ of type $A_1$ (see \cite{LaLu08}), and $1$ for $\Phi$ of
type $A_2$ (see Theorem~\ref{thm:A2}).  It remains a challenging
problem to find a conceptual interpretation of $\alpha_\Phi$ for
general $\Phi$.

For the proof of Theorem~\ref{thm:phi} and further details we refer to
the preprint~\cite{AvKlOnVo_arXiv}.



\end{article}
\end{document}